\documentclass[reqno,10pt]{amsart}
\usepackage{amsfonts}

\newtheorem{theorem}{Theorem}
\theoremstyle{plain}
\newtheorem{acknowledgement}{Acknowledgement}

\newtheorem{corollary}{Corollary}

\newtheorem{definition}{Definition}

\newtheorem{remark}{Remark}

\numberwithin{equation}{section}

\input{tcilatex}

\begin{document}
\author{}
\title{}
\maketitle

\begin{center}
\thispagestyle{empty}\textbf{{\Large {On }}}$q$\textbf{{\Large {-deformed
Stirling numbers}}}

\bigskip \bigskip

\textbf{{\Large {Yilmaz\ Simsek}}}

\textit{{\Large {University of Akdeniz, Faculty of Arts and Science,
Department of Mathematics,}}}{\Large {\textit{\ 07058 Antalya}, Turkey}}

{\Large 
}

{\Large \textbf{E-Mail: ysimsek@akdeniz.edu.tr,}}

{\Large 
}

\bigskip

\textbf{{\Large {Abstract}}}
\end{center}

The purpose of this article is to introduce $q$-deformed Stirling numbers of
the first and second kinds. Relations between these numbers, Riemann zeta
function and $q$-Bernoulli numbers of higher order are given. Some relations
related to the classical Stirling numbers and Bernoulli numbers of higher
order are found. By using derivative operator to the generating function of
the $q$-deformed Stirling numbers of the second kinds, a new function is
defined which interpolates the $q$-deformed Stirling numbers of the second
kinds at negative integers. The recurrence relations of the Stirling numbers
of the first and second kind are given. In addition, relation between $q$%
-deformed Stirling numbers and $q$-Bell numbers is obtained.\bigskip

\textbf{2000 Mathematics Subject Classification.} 11B39, 11B68,
11B73.\bigskip

\textbf{Key Words and Phrases.} $q$-Bernoulli numbers and polynomials, $q$%
-Stirling numbers first and second kind, fermionic Stirling numbers first
and second kind. $q$-Bell Numbers.\bigskip \bigskip

\section{Introduction, Definitions and Notations}

The $q$-deformed Stirling numbers of the first and second kind are denoted
by $s(n,k,q)$, $S(n,k,q)$, respectively. The fermionic Stirling numbers of
the first and second kind are denoted by $s_{f}(n,k)$, $S_{f}(n,k)$,
respectively. In this paper, we use notation in the work of Kim\cite{Kim5}
and Schork\cite{M. Schork}. $q$-Stirling numbers were first defined in the
work of Carlitz\cite{carlitz}. A lot of combinatorial work has centered
around the $q$-analogue, the earliest by Milne\cite{milne}; also see (\cite%
{M. Schork}, \cite{Ehrenborg}, \cite{Katriel}, \cite{Kim5}, \cite{kwaski}, 
\cite{kwasnki-1}, \cite{Wager}). In \cite{kimpjms}, Kim constructed $q$%
-Bernoulli numbers of higher order associated with the $p$-adic $q$%
-integers. $q$-Volkenborn integral is originately constructed by Kim\cite%
{Kim5}. By using the $q$-Volkenborn integral, Kim\cite{Kim5} evaluated
complete sum for $q$-Bernoulli polynomials. He also obtained relations
between $q$-Bernoulli numbers and $q$-analogs of the Stirling numbers. The
fermionic and bosonic Stirling numbers were given in detail by Schork\cite%
{M. Schork}, \cite{M. Schork-1}. In \cite{Wager}, Wagner studied three
partition statistics and the $q$-Stirling and $q$-Bell numbers that serve as
their generating functions, evaluating these numbers when $q=-1$.

In this paper, we give $q$-deformed Stirling numbers of the second kinds.
Relations between these numbers, Riemann zeta functions and $q$-Bernoulli
numbers higher order are given. We also give some relations related to the
classical Stirling numbers and Bernoulli numbers higher order. By using
derivative operator to the generating function of the $q$-deformed Stirling
numbers of the second kinds, we define new function, which interpolates the $%
q$-deformed Stirling numbers of the second kinds at negative integers. The
recurrence relations of the Stirling numbers of the first and the second
kind are given. We also give relation between $q$-deformed Stirling numbers
and $q$-Bell numbers.

Let $q\in (-1,1]$. The fundamental properties of the $q$-integers and $q$%
-deformed numbers are given by%
\begin{equation}
\lbrack n,q]=[n]=\frac{1-q^{n}}{1-q}.  \label{y1}
\end{equation}%
\begin{equation*}
\lbrack n]!=[n][n-1]...[2][1]\text{, }[0]!=1\text{ and }\left( 
\begin{array}{c}
n \\ 
k%
\end{array}%
\right) _{q}=\frac{[n]!}{[n-k]![k]!}.
\end{equation*}%
Note that $\lim_{q\rightarrow 1}[n]=n$, cf. (\cite{Katriel}, \cite%
{k1999jntqgamma}, \cite{k2002djm}, \cite{k2002padicinvar}, \cite{k2003nonars}%
, \cite{M. Schork}, \cite{k2004Sumsq}, \cite{tkim-skim-dwpark}, \cite{kwaski}%
).

The generating function of the $q$-Stirling numbers of the second kind is
given by defined by\cite{Kim5}:%
\begin{equation}
F_{S,q}(t)=\frac{q^{\frac{k(1-k)}{2}}}{[k]!}\sum_{j=0}^{k}(-1)^{k-j}\left( 
\begin{array}{c}
k \\ 
j%
\end{array}%
\right) _{q}q^{\frac{(k-j)(k-j-1)}{2}}e^{[j]t}=\sum_{n=0}^{\infty }\frac{%
S(n,k,q)t^{n}}{n!}.  \label{y2}
\end{equation}%
Let $(Eh)(x)=h(x+1)$ be the shift operator. Let $\Delta
_{q}^{n}=\prod_{j=0}^{n-1}(E-q^{j}I)$ be the $q$-difference operator cf. 
\cite{Kim5}.%
\begin{equation*}
\lbrack x]^{n}=\sum_{k=0}^{n}\left( 
\begin{array}{c}
x \\ 
k%
\end{array}%
\right) _{q}[k]!q^{\frac{k(1-k)}{2}}S(n,k,q)\text{ cf. (\cite{Kim5}),}
\end{equation*}%
where $S(n,k,q)$\ is denoted the $q$-Stirling numbers of the second kind. By
the above equation, we have%
\begin{equation*}
S(n,k,q)=q^{\frac{k(1-k)}{2}}\frac{1}{[k]!}\sum_{j=0}^{k}(-1)^{j}q^{\frac{%
j(j-1)}{2}}\left( 
\begin{array}{c}
k \\ 
j%
\end{array}%
\right) _{q}[k-j]^{n}\text{ cf. (\cite{Kim5}).}
\end{equation*}%
Kim\cite{Kim5} defined $q$-analog of the Newton-Gregory expansion as follows:%
\begin{equation*}
S(n,k,q)=\frac{q^{\frac{k(1-k)}{2}}}{[k]!}\Delta _{q}^{k}0^{n}
\end{equation*}

\section{$q$-deformed Stirling Numbers}

When $q<0$, we write $q\equiv -q^{\sim }$ with $q^{\sim }>0$. By (\ref{y1}),
we have%
\begin{equation*}
\lbrack n]=\frac{1-(-q^{\sim })^{n}}{1+q^{\sim }}=[n,q^{\sim }]^{F}\text{
cf. (\cite{M. Schork}),}
\end{equation*}%
$[n,q^{\sim }]^{F}$ is called $q^{\sim }$-fermionic basic numbers. These
numbers are appearing in recent studies of the $q^{\sim }$-deformed
fermionic oscillator (see \cite{kwasnki-1}, \cite{M. Schork}, \cite%
{Parthasarathy}, \cite{M. Schork-1}). Observe that $q\rightarrow -1$, i. e., 
$q^{\sim }\rightarrow 1$ yields (\cite{M. Schork})%
\begin{equation}
\lbrack n,q=-1]=[n,q^{\sim }=1]^{F}=\frac{1-(-1)^{n}}{2}=\epsilon
_{n}=\left\{ 
\begin{array}{c}
0\text{, if }n\text{ is even integers} \\ 
1\text{, if }n\text{ is odd integers.}%
\end{array}%
\right.  \label{ekk-3}
\end{equation}%
Hence, for $n\geq 2$ 
\begin{equation}
\lbrack n,q=-1]!=0.  \label{ek-3}
\end{equation}%
The $q^{\sim }$-fermionic basic numbers were given in detail by Schork\cite%
{M. Schork}, \cite{M. Schork-1}.

By using (\ref{y2}), the $q$-deformed Stirling numbers of the second kind, $%
S(n,k,q)$ are given by (in the version of Kim\cite{Kim5}):%
\begin{equation}
S(n,k,q)=\sum_{j=0}^{k}(-1)^{k-j}q^{\frac{k(1-k)+(k-j)(k-j-1)}{2}}\frac{%
[j]^{n-1}}{[j-1]![k-j]!}\text{.}  \label{ek-4}
\end{equation}%
where $n,k\in \mathbb{N}$ with $k\leq n$. The recurrence relations of $%
S(n,k,q)$, with $S(1,0,q)=0$ and $S(1,1,q)=1$, is given by%
\begin{equation*}
S(n+1,k,q)=q^{k-1}S(n,k-1,q)+[k]S(n,k,q)\text{, cf. (\cite{Charalambides},%
\cite{Ehrenborg},\cite{Kim5},\cite{milne}, \cite{M. Schork}).}
\end{equation*}

The recurrence relation implies that $S(n,k,q)$ are polynomials in $q$ (for
detail see \cite{M. Schork}). In the \textit{bosonic} limit: $%
\lim_{q\rightarrow 1}S(n,k,q)=S(n,k)$. Some special values of $k$, the $q$%
-deformed Stirling numbers of the second kind $S(n,k,q)$ are given by\cite%
{M. Schork}: $S(n,1,q)=1$, $S(n,2,q)=[2]^{n-1}-1$, $S(n,n,q)=q^{\frac{n(n-1)%
}{2}}$.

By applying the derivative operator $\frac{d^{k}}{dt^{k}}F_{S,q}(t)\mid
_{t=0}$ to (\ref{y2}), we arrive at the following theorem:

\begin{equation}
S(n,k,q)=\frac{d^{n}}{dt^{n}}F_{S,q}(t)\mid _{t=0}=\frac{q^{\frac{k(1-k)}{2}}%
}{[k]!}\sum_{j=0}^{k}(-1)^{k-j}\left( 
\begin{array}{c}
k \\ 
j%
\end{array}%
\right) _{q}q^{\frac{(k-j)(k-j-1)}{2}}[j]^{n}.  \label{a-1}
\end{equation}

By using (\ref{a-1}), we define $Y_{S}(z,k,q)$ function as follows:

\begin{definition}
Let $z\in \mathbb{C}$. We define%
\begin{equation}
Y_{S}(z,k,q)=\frac{q^{\frac{k(1-k)}{2}}}{[k]!}\sum_{j=0}^{k}(-1)^{k-j}\left( 
\begin{array}{c}
k \\ 
j%
\end{array}%
\right) _{q}q^{\frac{(k-j)(k-j-1)}{2}}[j]^{-z}.  \label{a-2}
\end{equation}
\end{definition}

Observe that if $z\in \mathbb{C}$, then $Y_{S}(z,k,q)$ is an analytic
function. The function $Y_{S}(z,k,q)$ interpolates the $q$-deformed Stirling
numbers of the second kind $S(n,k,q)$ at negative integers, which is given
in Theorem \ref{theoremA}, below.

By (\ref{y1}), (\ref{a-1}), we have the following Corollary:

\begin{corollary}
\begin{equation*}
S(n,k,q)=\frac{q^{\frac{k(1-k)}{2}}}{(1-q)^{n}[k]!}\sum_{j=0}^{k}%
\sum_{d=0}^{n}(-1)^{k-j-d}\left( 
\begin{array}{c}
k \\ 
j%
\end{array}%
\right) _{q}\left( 
\begin{array}{c}
n \\ 
d%
\end{array}%
\right) q^{\frac{(k-j)(k-j-1)+2jd}{2}}.
\end{equation*}
\end{corollary}

By substituting $z=-n$, with $n$ is a positive integer, into (\ref{a-2}),
and using (\ref{a-1}), we arrive at the following theorem:

\begin{theorem}
\label{theoremA}Let $n$ be a positive integer. Then we have%
\begin{equation}
Y_{S}(-n,k,q)=S(n,k,q).  \label{a-3}
\end{equation}
\end{theorem}

The $q$\textit{-deformed \textbf{Bell numbers} }are defined by \cite{M.
Schork}%
\begin{equation*}
B(n,q)=\sum_{k=0}^{n}S(n,k;q).
\end{equation*}

By using (\ref{a-1}) and (\ref{a-3}), we give relation between $Y_{S}(z,k,q)$
and $B(n,q)$ as follows:

\begin{theorem}
Let $n$ be a positive integer. Then we have%
\begin{eqnarray*}
B(n,q) &=&\sum_{k=0}^{n}Y_{S}(-n,k,q) \\
&=&\sum_{k=0}^{n}\sum_{j=0}^{k}(-1)^{k-j}\left( 
\begin{array}{c}
k \\ 
j%
\end{array}%
\right) _{q}\frac{[j]^{n}q^{\frac{k(1-k)+(k-j)(k-j-1)}{2}}}{[k]!}.
\end{eqnarray*}
\end{theorem}

\begin{remark}
Observe that when $\lim_{q\rightarrow 1}B(n,q)=B(n)=\sum_{k=0}^{n}S(n,k)$,
where $B(n)$ denotes the classical Bell numbers cf. (\cite{M. Schork}, \cite%
{milne}, \cite{Wager}). In \cite{Gessel}, Gessel gave relation between the
classical Stirling numbers of first kind, $S(n,k)$ and the classical
Bernoulli numbers of higher order,$B_{k}^{(n)}$ as follows:%
\begin{equation}
S(n+k,n)=\left( 
\begin{array}{c}
n+k \\ 
k%
\end{array}%
\right) B_{k}^{(-n)},  \label{a5}
\end{equation}%
where%
\begin{equation*}
\sum_{j=0}^{\infty }B_{j}^{(n)}\frac{t^{j}}{j!}=\left( \frac{t}{e^{t}-1}%
\right) ^{n}.
\end{equation*}%
In \cite{Kim5}, Kim gave relation between $S(n,k;q)$ numbers and $q$%
-Bernoulli numbers of higher-order as follows: Let $m\geq 0$ and $h,k$ are
natural numbers.%
\begin{equation*}
m^{k}\sum_{j=0}^{m}\left( 
\begin{array}{c}
m \\ 
j%
\end{array}%
\right) (q-1)^{j}\beta _{j}(o,-k,q)=\sum_{j=0}^{k}q^{\frac{k(k-1)}{2}%
}[j]!S(k,j;q)\left( 
\begin{array}{c}
m \\ 
j%
\end{array}%
\right) _{q},
\end{equation*}%
where%
\begin{equation*}
\beta _{m}(h,k,q)=(1-q)^{-m}\sum_{j=0}^{m}\left( 
\begin{array}{c}
m \\ 
j%
\end{array}%
\right) (-1)^{j}\left( \frac{h+j}{[h+j]}\right) ^{k}.
\end{equation*}%
In \cite{Rassias and Srivastava}, Rassias and Srivastava gave relation
between Riemann zeta functions and the classical Stirling numbers of first
kind, $s(n,k)$ as follows:%
\begin{equation}
\zeta (k+1)=\sum_{n=k}^{\infty }\frac{(-1)^{n-k}}{n.n!}s(n,k).  \label{a-4}
\end{equation}
\end{remark}

For each $k=0,1,...,n-1$, ($n\geq 1$), the \textit{Eulerian numbers} $E(n,k)$
are given by\cite{Cheon and Kim}%
\begin{equation*}
E(n,k)=\sum_{j=0}^{k}(-1)^{j}\left( 
\begin{array}{c}
n+1 \\ 
j%
\end{array}%
\right) (k+1-j)^{n}.
\end{equation*}%
Relation between $E(n,k)$ and the classical Stirling numbers of the second
kind, $S(n,k)$ is given by\cite{Cheon and Kim}%
\begin{equation}
S(n,m)=\frac{1}{m!}\sum_{j=0}^{n-1}E(n,j)\left( 
\begin{array}{c}
j \\ 
n-m%
\end{array}%
\right) ,\text{ }n\geq m,\text{ }n\geq 1.  \label{a-555}
\end{equation}

By (\ref{a5}) and (\ref{a-555}), after some elementary calculations, we
arrive at the following corollary:

\begin{corollary}
Let $k=0,1,...,n-1$, ($n\geq 1$). Then we have%
\begin{equation}
B_{k}^{(-n)}=\frac{\left( 
\begin{array}{c}
n+k \\ 
k%
\end{array}%
\right) }{n!}\sum_{j=0}^{n+k-1}E(n+k,j)\left( 
\begin{array}{c}
j \\ 
k%
\end{array}%
\right) .  \label{a-6}
\end{equation}
\end{corollary}

By (\ref{a-3}) with $q\rightarrow 1$ and (\ref{a-555}), we have the
following Corollary:

\begin{corollary}
Let $n$ be a positive integer and $n\geq k,$ $n\geq 1$. Then we have%
\begin{equation*}
Y_{S}(-n,k)=\frac{1}{k!}\sum_{j=0}^{n-1}E(n,j)\left( 
\begin{array}{c}
j \\ 
n-k%
\end{array}%
\right) .
\end{equation*}
\end{corollary}

The $q$-deformed Stirling number of the first kind $s(n,k,q)$, with $%
s(1,0,q)=0$ and $s(1,1,q)=1$, satisfy the following recurrence relations (
see \cite{M. Schork}, \cite{Charalambides}):%
\begin{equation}
s(n+1,k,q)=q^{-n}\left( s(n,k-1,q)-[n]s(n,k,q)\right) .  \label{ek-7}
\end{equation}%
If $q\rightarrow 1$ in the above,\ we have%
\begin{equation*}
s(n+1,k)=s(n,k-1)-ns(n,k)
\end{equation*}

For $k\geq 1$, we set%
\begin{equation*}
\lbrack x]^{\underline{k}}=[x][x-1]...[x-k+1],\text{ cf. (\cite{carlitz}, 
\cite{Katriel}, \cite{M. Schork}, \cite{M. Schork-1}, \cite{Parthasarathy}).}
\end{equation*}%
By using the above relation, the $q$-deformed Stirling numbers is defined as 
\cite{M. Schork}%
\begin{equation}
\lbrack x]^{n}=\sum_{j=0}^{n}S(n,j,q)[x]^{\underline{j}}\text{ and }[x]^{%
\underline{n}}=\sum_{j=0}^{n}s(n,j,q)[x]^{j}.  \label{ek-8}
\end{equation}%
By using (\ref{ek-8}), the $q$-deformed Stirling numbers of the first and
the second kind satisfies for $n\geq m$ the inversion relations, which are
given by the following theorem:

\begin{theorem}
Let $n$ and $m$ be non-negative integers. Then we have%
\begin{eqnarray*}
\sum_{k=m}^{n}s(n,k,q)S(k,m,q) &=&\delta _{n,m}, \\
\sum_{k=m}^{n}S(n,k,q)s(k,m,q) &=&\delta _{n,m}.
\end{eqnarray*}
\end{theorem}

Proofs of this theorem were given by Schork\cite{M. Schork} and Charalambides%
\cite{Charalambides}. From the above theorem, $q$-Stirling numbers of the
first and the second kind satisfy the ortogonality relations.

\section{Further Remarks and Observations}

The fermionic Stirling numbers of the first and the second kind studied by 
\cite{milne}, \cite{M. Schork}, \cite{M. Schork-1}, \cite{carlitz}, \cite%
{Ehrenborg}, \cite{Wager}, \cite{kwaski}, \cite{kwasnki-1}, \cite%
{Parthasarathy}, \cite{Katriel}, \cite{Kim5}. In this section, we can use
some notations which are due to Schork\cite{M. Schork}, and Kim\cite{Kim5}.
In \cite{M. Schork}, Schork gave the recurrence relations of the fermionic
Stirling numbers of first kind and second kind as follows:%
\begin{equation}
s_{f}(n+1,k)=(-1)^{n}s_{f}(n,k-1)+(-1)^{n+1}\epsilon _{n}s_{f}(n,k),
\label{ek-10}
\end{equation}%
with $s_{f}(1,0)=0$ and $s_{f}(1,1)=1$. For the convention, here we take $%
s_{f}(n,0)=0$ and%
\begin{equation}
S_{f}(n+1,k)=(-1)^{k-1}S_{f}(n,k-1)+\epsilon _{k}S_{f}(n,k),  \label{ek-11}
\end{equation}%
with $S_{f}(1,0)=0$ and $S_{f}(1,1)=1$, where $\epsilon _{k}$ is defined in (%
\ref{ekk-3}). Schork gave the values $S_{f}(n,k)$ for maximal and small $k$
and found $S_{f}(n,n)=(-1)^{\frac{n(n-1)}{2}}$ as well as $S_{f}(n,1)=1$, $%
S_{f}(n,2)=-1$, $S_{f}(n,3)=2-n$, $S_{f}(n,4)=n-3$ and by (\ref{ek-10}), we
easily see that $s_{f}(n,n)=(-1)^{\frac{n(n-1)}{2}}$.

By (\ref{ek-10}) and (\ref{ekk-3}), many Stirling numbers of the first kind
vanish. By induction over $k$, Schork\cite{M. Schork} prove $s_{f}(n,k)=0$
for $n>2k$. By the same method, $s_{f}(n,k^{\sim })=0$ for $1\leq k^{\sim
}\leq k$ and $n>2k^{\sim }$ and $s_{f}(n,k+1)=0$ for $n\geq 2k+3$ (for
detail see \cite{M. Schork}). From (\ref{ek-10}), we arrive at the following
theorem\cite{M. Schork}:

\begin{theorem}
Let $k$ be non-negative integer. Then we have%
\begin{equation*}
s_{f}(2k+3,k+1)=s_{f}(2k+2,k)-\epsilon _{2k+2}s_{f}(2k+2,k+1).
\end{equation*}
\end{theorem}

Note that by the induction hypothesis, the first summand vanishes, whereas
by the $\epsilon _{2k+2}=0$, the second summand vanish. Hence, for a given $%
n $, the first $\frac{n}{2}$ Stirling numbers $s_{f}(n,k)$ vanish. Let%
\begin{equation*}
\lbrack n]_{f}=[n,q=-1]=\epsilon _{n}\text{ cf. \cite{M. Schork}}
\end{equation*}%
for the fermionic basic numbers, by (\ref{ekk-3}). The fermionic Stirling
numbers are connection coefficients for the fermionic basic numbers \cite{M.
Schork}:%
\begin{equation}
\lbrack x]_{f}^{n}=\sum_{j=0}^{n}S_{f}(n,j)[x]_{f}^{\underline{j}}\text{ and 
}[x]_{f}^{\underline{n}}=\sum_{j=0}^{n}s_{f}(n,j)[x]_{f}^{j}.  \label{ek-12}
\end{equation}%
By induction over $n$, Schork\cite{M. Schork} proved the first equation. We
give sketch of the proof as follows: Let the assertion obtain for $n$. Thus,
the induction hypothesis implies that%
\begin{equation*}
\lbrack x]_{f}^{n+1}=\left( \sum_{j=0}^{n}S_{f}(n,j)[x]_{f}^{\underline{j}%
}\right) [x]_{f}.
\end{equation*}%
By using $f$-arithmetic relations, $[x]_{f}=[j]_{f}+(-1)^{j}[x-j]_{f}$, $%
[n+m]_{f}=[m]_{f}+(-1)^{m}[n]_{f}$ and $[n+1]_{f}=1-[n]_{f}$,\ we easily find%
\begin{eqnarray}
\lbrack x]_{f}^{n+1} &=&\sum_{j=0}^{n}S_{f}(n,j)[x]_{f}^{\underline{j}%
}[x]_{f}=\sum_{j=0}^{n}\left(
S_{f}(n,j)[j]_{f}+(-1)^{j}[x-k]_{f}S_{f}(n,j)\right) [x]_{f}^{\underline{j}}
\notag \\
&=&\sum_{j=0}^{n}S_{f}(n,j)[j]_{f}[x]_{f}^{\underline{j}}+%
\sum_{j=0}^{n}(-1)^{j}[x-k]_{f}S_{f}(n,j)[x]_{f}^{\underline{j}}  \notag \\
&=&\sum_{j=0}^{n}S_{f}(n,j)[j]_{f}[x]_{f}^{\underline{j}}+%
\sum_{j=0}^{n}(-1)^{j}S_{f}(n,j)[x]_{f}^{\underline{j+1}}  \notag \\
&=&\sum_{j=1}^{n+1}\left( S_{f}(n,j)[j]_{f}+(-1)^{j-1}S_{f}(n,j-1)\right)
[x]_{f}^{\underline{j}}.  \label{ek-a}
\end{eqnarray}%
Since $S_{f}(n+1,0)=0$, we have%
\begin{equation}
\lbrack x]_{f}^{n+1}=\sum_{j=0}^{n+1}S_{f}(n+1,j)[x]_{f}^{\underline{j}%
}=\sum_{j=1}^{n+1}S_{f}(n+1,j)[x]_{f}^{\underline{j}}.  \label{ek-b}
\end{equation}%
By (\ref{ek-a}) and (\ref{ek-b}), we arrive at (\ref{ek-11}).

Since $[x-n]_{f}=(-1)^{n}[x]_{f}+(-1)^{n+1}\epsilon _{n}$, we have

\begin{eqnarray*}
\lbrack x]_{f}^{\underline{n+1}}
&=&[x-n]_{f}\sum_{j=0}^{n}s_{f}(n,j)[x]_{f}^{j}=\sum_{j=0}^{n}s_{f}(n,j)%
\left( (-1)^{n}[x]_{f}+(-1)^{n+1}\epsilon _{n}\right) [x]_{f}^{j} \\
&=&\sum_{j=0}^{n}(-1)^{n}s_{f}(n,j)[x]_{f}^{j+1}+%
\sum_{j=0}^{n}(-1)^{n+1}s_{f}(n,j)\epsilon _{n}[x]_{f}^{j} \\
&=&\sum_{j=1}^{n+1}\left( (-1)^{n}s_{f}(n,j-1)+(-1)^{n+1}s_{f}(n,j)\epsilon
_{n}\right) [x]_{f}^{j},
\end{eqnarray*}%
by the above equation, we obtain (\ref{ek-10}).

Observe that by (\ref{ek-10}) and (\ref{ek-11}), Schork\cite{M. Schork}
defined the fermionic Stirling numbers of first and second kind,
respectively.

It is well-known that if $j\geq 2$, then $[x]_{f}^{\underline{j}}$ vanishes,
so by using the first equation of (\ref{ek-12})%
\begin{equation*}
\lbrack x]_{f}^{n}=S_{f}(n,1)[x]_{f}=[x]_{f}.
\end{equation*}

It is easy to prove the following theorem.

\begin{theorem}
Let $n$ and $m$ be non-negative integers. Then we have%
\begin{equation}
\sum_{j=m}^{n}s_{f}(n,j)S_{f}(j,m)=\delta _{n,m},\text{ and }%
\sum_{j=m}^{n}S_{f}(n,j)s_{f}(j,m)=\delta _{n,m}.  \label{ek-13}
\end{equation}
\end{theorem}

By (\ref{ek-13}), the fermionic Stirling numbers satisfy the inversion
relations. The \ proof of (\ref{ek-13}), was given by Schork\cite{M. Schork}%
. He also proved the following the recurrence relations of the fermionic
Stirling numbers of second kind: For odd $j>3$ one can find that $%
S_{f}(n+1,j)=S_{f}(n,j)+S_{f}(n,j-1)$. Let $n$ and $j$ be non-negative
integers. 
\begin{equation*}
S_{f}(n+1,j)=S_{f}(n,j)-S_{f}(n-1,j-2).
\end{equation*}

\begin{acknowledgement}
We would like to thank to Professor Matthias Schork for his many valuable
suggestions and comments on this paper. We also would like to thank to
Professor Taekyun Kim for his many valuable comments on this paper.

This work was supported by Akdeniz University Scientific Research Projects
Unit.
\end{acknowledgement}

\end{document}